\documentclass{amsproc}%
\usepackage{amsfonts}
\usepackage{amsmath}
\usepackage{amssymb}
\usepackage{graphicx}%
\theoremstyle{plain}

\newtheorem{algorithm}{Algorithm}

\newtheorem{definition}{Definition}

\newtheorem{lemma}{Lemma}

\newtheorem{proposition}{Proposition}

\newtheorem{theorem}{Theorem}
\numberwithin{equation}{section}
\begin{document}
\title[Hook-lengths]{Hook-lengths and Pairs of Compositions}
\author{Charles F. Dunkl}
\address{Department of Mathematics, P.O. Box 400137\\
University of Virginia\\
Charlottesville, VA 22904-4137 U.S.}
\email{cfd5z@virginia.edu}
\urladdr{http://www.people.virginia.edu/\symbol{126}cfd5z}
\thanks{During the preparation of this article the author was partially supported by
NSF grant DMS 0100539}
\date{May 20, 2005}

\begin{abstract}
The monomial basis for polynomials in N variables is labeled by compositions.
To each composition there is associated a hook-length product, which is a
product of linear functions of a parameter. The zeroes of this product are
related to "critical pairs" of compositions; a concept defined in this paper.
This property can be described in an elementary geometric way; for example:
consider the two compositions (2,7,8,2,0,0) and (5,1,2,5,3,3), then the
respective ranks, permutations of the index set \{1,2,...,6\} sorting the
compositions, are (3,2,1,4,5,6) and (1,6,5,2,3,4), and the two vectors of
differences (between the compositions and the ranks, respectively) are
(-3,6,6,-3,-3,-3) and (2,-4,-4,2,2,2), which are parallel, with ratio -3/2.
For a given composition and zero of its hook-length product there is an
algorithm for constructing another composition with the parallelism property
and which is comparable to it in a certain partial order on compositions,
derived from the dominance order. This paper presents the motivation from the
theory of nonsymmetric Jack polynomials and the description of the algorithm,
as well as the proof of its validity.

\end{abstract}
\maketitle

\section{Introduction}

A composition is an element of $\mathbb{N}_{0}^{N}$ (where $\mathbb{N}%
_{0}:=\left\{  0,1,2,3,\ldots\right\}  $); a typical composition is
$\alpha=\left(  \alpha_{1},\ldots,\alpha_{N}\right)  $ and the components
$\alpha_{i}$ are called the parts of $\alpha$. Compositions have the obvious
application of labeling the monomial basis of polynomials in the variables
$x_{1},\ldots,x_{N}$ and they also serve as labels for the nonsymmetric Jack
polynomials (a set of homogeneous polynomials which are simultaneous
eigenfunctions of a certain parametrized and commuting set $\left\{
\mathcal{U}_{i}:1\leq i\leq N\right\}  $ of difference-differential
operators). In this context the ranks of the parts of a composition become
significant. The ranks are based on sorting on magnitude and index so that the
largest part has rank 1; if a value is repeated then the one with lower index
has lower rank. This is made precise in the following (the cardinality of a
set $E$ is denoted by $\#E$):

\begin{definition}
For $\alpha\in\mathbb{N}_{0}^{N}$ and $1\leq i\leq N$ let $r\left(
\alpha,i\right)  :=\#\left\{  j:\alpha_{j}>\alpha_{i}\right\}  +\#\left\{
j:1\leq j\leq i,\alpha_{j}=\alpha_{i}\right\}  $ be the rank function.
\end{definition}

A consequence of the definition is that $r\left(  \alpha,i\right)  <r\left(
\alpha,j\right)  $ is equivalent to $\alpha_{i}>\alpha_{j}$, or $\alpha
_{i}=\alpha_{j}$ and $i<j$. For any $\alpha$ the function $i\mapsto r\left(
\alpha,i\right)  $ is one-to-one on $\left\{  1,2,\ldots,N\right\}  $. A
partition is a composition satisfying $\alpha_{i}\geq\alpha_{i+1}$ for all
$i$, equivalently, $r\left(  \alpha,i\right)  =i$ for all $i$. For a fixed
$\alpha\in\mathbb{N}_{0}^{N}$ the values $\left\{  r\left(  \alpha,i\right)
:1\leq i\leq N\right\}  $ are independent of trailing zeros, that is, if
$\alpha^{\prime}\in\mathbb{N}_{0}^{M},\alpha_{i}^{\prime}=\alpha_{i}$ for
$1\leq i\leq N$ and $\alpha_{i}^{\prime}=0$ for $N<i\leq M$ then $r\left(
\alpha,i\right)  =r\left(  \alpha^{\prime},i\right)  $ for $1\leq i\leq N$,
and $r\left(  \alpha^{\prime},i\right)  =i$ for $N<i\leq M$. A formal
parameter $\kappa$ appears in the construction of nonsymmetric Jack
polynomials; their coefficients are in $\mathbb{Q}\left(  \kappa\right)  $, a
transcendental extension of $\mathbb{Q}$. The relevant information in a
composition label is encoded as the function $i\mapsto\alpha_{i}-\kappa
r\left(  \alpha,i\right)  $. We will be concerned with situations where a pair
$\left(  \alpha,\beta\right)  $ of compositions has the property that
$\alpha_{i}-\kappa r\left(  \alpha,i\right)  =\beta_{i}-\kappa r\left(
\beta,i\right)  $ for all $i$, when $\kappa$ is specialized to some negative
rational number. This is equivalent to the condition that $\left(  r\left(
\beta,i\right)  -r\left(  \alpha,i\right)  \right)  \kappa+\alpha_{i}%
-\beta_{i}$ is a rational multiple of $m\kappa+n$ for some fixed $m,n>0$ (or
that the vectors $\left(  \alpha_{i}-\beta_{i}\right)  _{i=1}^{N}$ and
$\left(  r\left(  \alpha,i\right)  -r\left(  \beta,i\right)  \right)
_{i=1}^{N}$ are parallel). For our application an additional condition is
imposed on the pair $\left(  \alpha,\beta\right)  $ which is stated in terms
of a partial order on compositions. Let $S_{N}$ denote the symmetric group on
$N$ objects, considered as the permutation group of $\left\{  1,2,\ldots
,N\right\}  $. The action of $S_{N}$ on compositions is defined by $\left(
w\alpha\right)  _{i}=\alpha_{w^{-1}\left(  i\right)  },1\leq i\leq N$.

\begin{definition}
For a composition $\alpha\in\mathbb{N}_{0}^{N}$ let $\left\vert \alpha
\right\vert :=\sum_{i=1}^{N}\alpha_{i}$ and let $\ell\left(  \alpha\right)
:=\max\left\{  j:\alpha_{j}>0\right\}  $ be the length of $\alpha$.
\end{definition}

\begin{definition}
For $\alpha\in\mathbb{N}_{0}^{N}$ let $\alpha^{+}$ denote the unique partition
such that $\alpha^{+}=w\alpha$ for some $w\in S_{N}$. For $\alpha,\beta
\in\mathbb{N}_{0}^{N}$ the partial order $\alpha\succ\beta$ ($\alpha$
dominates $\beta$) means that $\alpha\neq\beta$ and $\sum_{i=1}^{j}\alpha
_{i}\geq\sum_{i=1}^{j}\beta_{i}$ for $1\leq j\leq N$; and $\alpha
\vartriangleright\beta$ means that $\left\vert \alpha\right\vert =\left\vert
\beta\right\vert $ and either $\alpha^{+}\succ\beta^{+}$ or $\alpha^{+}%
=\beta^{+}$ and $\alpha\succ\beta$.
\end{definition}

For a given $\alpha\in\mathbb{N}_{0}^{N}$ let $w$ be the inverse function of
$i\mapsto r\left(  \alpha,i\right)  $ then $r\left(  \alpha,w\left(  j\right)
\right)  =j$ for $1\leq j\leq N$ and $\alpha=w\alpha^{+}$. This permutation
appears again in part (iii) of Proposition \ref{a-tilde}.

\begin{definition}
A pair $\left(  \alpha,\beta\right)  $ of compositions is a $\left(  -\frac
{n}{m}\right)  $-critical pair (where $m,n\geq1$) if $\alpha\vartriangleright
\beta$ and $m\kappa+n$ divides $\left(  r\left(  \beta,i\right)  -r\left(
\alpha,i\right)  \right)  \kappa+\alpha_{i}-\beta_{i}$ \ (in $\mathbb{Q}%
\left[  \kappa\right]  $) for each $i$.
\end{definition}

The divisibility property is equivalent to $\left(  r\left(  \beta,i\right)
-r\left(  \alpha,i\right)  \right)  n=\allowbreak m\left(  \alpha_{i}%
-\beta_{i}\right)  $ for all $i$. By elementary arguments we show why only
negative numbers appear in the critical pairs, and we also find a bound on
$\ell\left(  \beta\right)  $. A simple example shows that $m=0$ is possible:
let $\alpha=\left(  3,0\right)  $ and $\beta=\left(  2,1\right)  $, then both
$\alpha$ and $\beta$ have ranks $\left(  1,2\right)  .$

\begin{proposition}
Suppose $\alpha,\beta\in\mathbb{N}_{0}^{N},\alpha\vartriangleright\beta$ and
there are integers $m,n$ such that $\left(  \left(  r\left(  \beta,i\right)
-r\left(  \alpha,i\right)  \right)  \kappa+\alpha_{i}-\beta_{i}\right)
/\left(  m\kappa+n\right)  \in\mathbb{Q}$ for $1\leq i\leq N$, then $mn\geq0$
and $n\neq0$.
\end{proposition}

\begin{proof}
The case $n=0$ is impossible since that would imply $\alpha_{i}-\beta_{i}=0$
for all $i$, that is, $\alpha=\beta$. So we assume $n\geq1$ and then show
$m\geq0$. Let $w$ be the inverse function of $i\mapsto r\left(  \alpha
,i\right)  $ (so that $r\left(  \alpha,w\left(  i\right)  \right)  =i$). By
definition either $\alpha^{+}\succ\beta^{+}$ or $\alpha^{+}=\beta^{+}$ and
$\alpha\succ\beta$. Suppose that $\alpha^{+}\succ\beta^{+}$ and let $k\geq1$
have the property that $\beta_{w\left(  j\right)  }=\alpha_{w\left(  j\right)
}$ and $r\left(  \beta,w\left(  j\right)  \right)  =j$ for $1\leq j<k$ and at
least one of $\beta_{w\left(  k\right)  }\neq\alpha_{w\left(  k\right)  }$ and
$r\left(  \beta,w\left(  k\right)  \right)  >k=r\left(  \alpha,w\left(
k\right)  \right)  $ holds. Define $l$ by $r\left(  \beta,l\right)  =k$, then
by the definition of the dominance order $\succ$ we have that $\alpha
_{w\left(  k\right)  }\geq\beta_{l}$. Also $\beta_{l}\geq\beta_{w\left(
k\right)  }$ because $r\left(  \beta,w\left(  k\right)  \right)  \geq k$. The
case $\alpha_{w\left(  k\right)  }=\beta_{w\left(  k\right)  }$ and $r\left(
\beta,k\right)  >k$ (thus $n=0$) is impossible hence $\alpha_{w\left(
k\right)  }>\beta_{w\left(  k\right)  }$. If $r\left(  \beta,k\right)  =k$
then $m=0$ or else $r\left(  \beta,k\right)  >k$ and $m>0$.

Now suppose $\alpha^{+}=\beta^{+}$ and $\alpha\succ\beta$, and let $k\geq1$
have the property that $\beta_{j}=\alpha_{j}$ for $1\leq j<k$ and $\alpha
_{k}>\beta_{k}$ (the existence of $k$ follows from the definition of
$\alpha\succ\beta$). Since $\beta$ is a permutation of $\alpha$ we have that
$r\left(  a,j\right)  =r\left(  \beta,j\right)  $ for $1\leq j<k$ and
$r\left(  \alpha,k\right)  <r\left(  \beta,k\right)  $. This implies $m>0$.
\end{proof}

\begin{proposition}
\label{lgt}Suppose that $\left(  \alpha,\beta\right)  $ is a $\left(
-\frac{n}{m}\right)  $-critical pair, for some $m,n\geq1$, then $\ell\left(
\beta\right)  \leq\ell\left(  \alpha\right)  +\left\vert \alpha\right\vert $.
\end{proposition}

\begin{proof}
First we show that if $i>\ell\left(  \alpha\right)  $ and $\beta_{i}=0$ then
$\beta_{j}=0$ for all $j>i$. By hypothesis $\left(  m\kappa+n\right)  $
divides $\left(  r\left(  \beta,i\right)  -r\left(  \alpha,i\right)  \right)
\kappa+(\alpha_{i}-\beta_{i})=\left(  r\left(  \beta,i\right)  -i\right)
\kappa$, hence $r\left(  \beta,i\right)  =i$. This implies that $0\leq
\beta_{j}\leq\beta_{i}=0$ for all $j>i$. Thus if $\ell\left(  \beta\right)
>\ell\left(  \alpha\right)  $ then $\beta_{i}\geq1$ for $\ell\left(
\alpha\right)  <i\leq\ell\left(  \beta\right)  $. Since $\left\vert
\beta\right\vert =\left\vert \alpha\right\vert $ this shows that $\ell\left(
\beta\right)  -\ell\left(  \alpha\right)  \leq\left\vert \alpha\right\vert $.
\end{proof}

The motivation for the concept of hook-lengths associated with a composition
came from the representation theory of the symmetric group, where it appeared
in the famous hook-length formula for the degree of an irreducible
representation. In Section 2 we will explain the connection with nonsymmetric
Jack polynomials. However the following definitions are logically independent
of this theory. Suppose $\alpha\in\mathbb{N}_{0}^{N}$ and $\ell\left(
\alpha\right)  =m$; the (modified for compositions) \textit{Ferrers diagram}
of $\alpha$ is the set $\left\{  \left(  i,j\right)  :1\leq i\leq m,0\leq
j\leq\alpha_{i}\right\}  .$ For each node $\left(  i,j\right)  $ with $1\leq
j\leq\alpha_{i}$ \ there are two special subsets of the Ferrers diagram, the
\textit{arm} $\left\{  \left(  i,l\right)  :j<l\leq\alpha_{i}\right\}  $ and
the \textit{leg} $\left\{  \left(  l,j\right)  :l>i,j\leq\alpha_{l}\leq
\alpha_{i}\right\}  \cup\left\{  \left(  l,j-1\right)  :l<i,j-1\leq\alpha
_{l}<\alpha_{i}\right\}  $. The node itself, the arm and the leg make up the
\textit{hook}. (Note that for the case of partitions the nodes $\left(
i,0\right)  $ are omitted from the Ferrers diagram.)

Here is an example: the Ferrers diagram for the composition $\alpha=\left(
1,0,5,3,4,2\right)  $ (where the first part corresponds to the top row) is
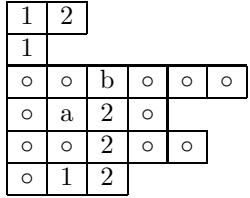
\[%
\begin{tabular}
[c]{|c|ccccc}\cline{1-2}%
1 & 2 & \multicolumn{1}{|c}{} &  &  & \\\cline{1-2}%
1 &  &  &  &  & \\\hline
$\circ$ & $\circ$ & \multicolumn{1}{|c}{b} & \multicolumn{1}{|c}{$\circ$} &
\multicolumn{1}{|c}{$\circ$} & \multicolumn{1}{|c|}{$\circ$}\\\hline
$\circ$ & a & \multicolumn{1}{|c}{2} & \multicolumn{1}{|c}{$\circ$} &
\multicolumn{1}{|c}{} & \\\cline{1-5}%
$\circ$ & $\circ$ & \multicolumn{1}{|c}{2} & \multicolumn{1}{|c}{$\circ$} &
\multicolumn{1}{|c}{$\circ$} & \multicolumn{1}{|c}{}\\\cline{1-5}%
$\circ$ & 1 & \multicolumn{1}{|c}{2} & \multicolumn{1}{|c}{} &  &
\\\cline{1-3}%
\end{tabular}
\]
The leg of the node $\left(  4,1\right)  $, labeled \textquotedblleft
a\textquotedblright, consists of the nodes labeled \textquotedblleft%
1\textquotedblright, and the leg of the node $\left(  3,2\right)  $, labeled
\textquotedblleft b\textquotedblright, consists of the nodes labeled
\textquotedblleft2\textquotedblright.

The cardinality of the leg is called the leg-length, formalized by the following:

\begin{definition}
For $\alpha\in\mathbb{N}_{0}^{N},1\leq i\leq\ell\left(  \alpha\right)  $ and
$1\leq j\leq\alpha_{i}$ the leg-length is%
\[
L\left(  \alpha;i,j\right)  :=\#\left\{  l:l>i,j\leq\alpha_{l}\leq\alpha
_{i}\right\}  +\#\left\{  l:l<i,j\leq\alpha_{l}+1\leq\alpha_{i}\right\}  .
\]

\end{definition}

For $t\in\mathbb{Q}\left(  \kappa\right)  $ the \textit{hook-length }and the
hook-length product for $\alpha$ are given by
\begin{align*}
h\left(  \alpha,t;i,j\right)   &  :=\left(  \alpha_{i}-j+t+\kappa L\left(
\alpha;i,j\right)  \right) \\
h\left(  \alpha,t\right)   &  :=\prod_{i=1}^{\ell\left(  \alpha\right)  }%
\prod_{j=1}^{\alpha_{i}}h\left(  \alpha,t;i,j\right)  ,
\end{align*}

Note that the indices $\left\{  i:\alpha_{i}=0\right\}  $ are omitted in the
product $h\left(  \alpha,t\right)  $. (In the present paper $t$ almost always
has the value $\kappa+1$, but $t=1$ and $t=\kappa$ do occur in some formulae
for Jack polynomials.) The results of Knop and Sahi \cite{KS} imply that for
any node $\left(  i,j\right)  ,1\leq j\leq\alpha_{i}$, so that $h\left(
\alpha,\kappa+1;i,j\right)  =\left(  L\left(  \alpha;i,j\right)  +1\right)
\kappa+\alpha_{i}+1-j$, there must exist at least one $\beta$ such that
$\left(  \alpha,\beta\right)  $ is $-\left(  \alpha_{i}+1-j\right)  /\left(
L\left(  \alpha;i,j\right)  +1\right)  $-critical. The main purpose of this
paper is to construct such a composition $\beta$ by direct algorithmic means.
This forms the content of Section 3. There are examples and discussion of open
problems in Section 4.

First we assume that the node is in the largest part, that is, $r\left(
\alpha,i\right)  =1$. The modification for other parts is trivial - one merely
ignores all larger parts $\alpha_{k}$ (with $r\left(  \alpha,k\right)
<r\left(  \alpha,i\right)  $). This will be explained in detail later. We
illustrate how the algorithm works on a partition $\alpha$, thereby avoiding
some technical complexity. Choose $n\leq\alpha_{1}$ and suppose $L\left(
\alpha;1,\alpha_{1}+1-n\right)  =m-1$ \ (thus $h\left(  \alpha,\kappa
+1;1,\alpha_{1}+1-n\right)  =m\kappa+n$). Then $\alpha_{m}>\alpha_{1}%
-n\geq\alpha_{m+1}$. Define a sequence by $\xi_{mk+i}=\alpha_{i}-nk$ for
$1\leq i\leq m$ and $k\geq0$. Then $\xi_{j}\geq\xi_{j+1}$ for all $j$; indeed
if $j=mk+i$ with $i<m$ then $\xi_{j}-\xi_{j+1}=\left(  \alpha_{i}-nk\right)
-\left(  \alpha_{i+1}-nk\right)  \geq0$ and if $j=mk$ then $\xi_{mk}%
-\xi_{mk+1}=\left(  \alpha_{m}-n\left(  k-1\right)  \right)  -\left(
\alpha_{1}-nk\right)  =\alpha_{m}-\left(  \alpha_{1}-n\right)  >0$. Also
$\alpha_{m+1}\leq\alpha_{1}-n=\xi_{m+1}$. Since the values $\left\{  \xi
_{j}\right\}  $ are eventually negative there exists a unique $T$ such that
$\alpha_{m+s}\leq\xi_{m+s+1}$ for $0<s<T$ and $\alpha_{m+T}>\xi_{m+T+1}$ (or
$T=1$ when $\alpha_{m+1}>\xi_{m+2}$). Set $t:=\left(  \left(  T-1\right)
\operatorname{mod}m\right)  +1$ and $k:=\left(  T-t\right)  /m$ (thus $T=mk+t$
and $1\leq t\leq m$). Define $\beta\in\mathbb{N}_{0}^{N}$ by
\[
\beta_{i}:=\left\{
\begin{tabular}
[c]{ll}%
$\xi_{\left(  k+1\right)  m+i}=\alpha_{i}-\left(  k+1\right)  n,$ & $1\leq
i\leq t$\\
$\xi_{km+i}=\alpha_{i}-kn,$ & $t<i\leq m$\\
$\alpha_{i}+n,$ & $m+1\leq i\leq m+T$\\
$\alpha_{i},$ & $m+T<i.$%
\end{tabular}
\ \ \ \ \ \right.
\]
In this context, an upper bound on $N$ is not needed; that is, $\alpha_{i}$ is
defined for all $i\geq1$ and $\alpha_{i}=0$ for $i>\ell\left(  \alpha\right)
$. However one can show that $l\left(  \beta\right)  \leq\max\left(  l\left(
\alpha\right)  ,m+T_{0}\right)  $ where $T_{0}:=\sum_{i=1}^{m}\left\lfloor
\frac{\alpha_{i}}{n}\right\rfloor $ ($\left\lfloor r\right\rfloor $ is the
largest integer $\leq r$) and it suffices to take $N$ as large as this bound.
Then $\alpha_{m+T}+n=\beta_{m+T}>\beta_{t+1}\geq\ldots\geq\beta_{m}>\beta
_{1}\geq\ldots\geq\beta_{t}\geq\beta_{m+T+1}=\alpha_{m+T+1}$ and $r\left(
\beta,m+i\right)  =i$ for $1\leq i\leq T$, $r\left(  \beta,i\right)
=T+m-t+i=m\left(  k+1\right)  +i$ for $1\leq i\leq t,r\left(  \beta,i\right)
=mk+i$ for $t<i\leq m$ and $r\left(  \beta,i\right)  =i$ for $i>m+T$. The
proof of these facts is a special case of the general result.

The computational scheme can be set up in algorithmic fashion:\ consider the
example $\alpha=\left(  9,8,8,7,4,3,3,2,2\right)  $ with $h\left(
\alpha,\kappa+1;1,7\right)  =4\kappa+3$, so $m=4$ and $n=3$. Generate enough
of the sequence $\left(  \xi_{i}\right)  _{i=1}^{N}$ to determine the value of
$T$; note that $\left(  \xi_{i}\right)  _{i=5}^{N}=\left(
9-3,8-3,8-3,7-3,9-6,8-6,\ldots\right)  $. Comparing the sequences%
\begin{align*}
\alpha &  =\left(  9,8,8,7,4,3,3,2,2,0,0,0,0,0,\ldots\right) \\
\xi &  =\left(  9,8,8,7,6,5,5,4,3,2,2,1,0,-1,\ldots\right)
\end{align*}
term-by-term we see that $T=9$ ($\alpha_{4+9}>\xi_{4+10}$ and $\alpha
_{4+s}\leq\xi_{5+s}$ for $1\leq s\leq8$.) Finally $t=1,k=2$ and the formula
produces $\beta=\left(  0,2,2,1,7,6,6,5,5,3,3,3,3\right)  $; it can be checked
that $\left(  \alpha,\beta\right)  $ is a $\left(  -\frac{3}{4}\right)
$-critical pair.

\section{Nonsymmetric Jack polynomials and hook-length products}

For $\alpha\in\mathbb{N}_{0}^{N}$ the corresponding monomial is $x^{\alpha
}:=\prod_{i=1}^{N}x_{i}^{\alpha_{i}}$ and the degree of $x^{\alpha}$ is
$\left\vert \alpha\right\vert $. For $1\leq i,j\leq N$ and $i\neq j$ the
transposition of $i$ and $j$ is denoted by $\left(  i,j\right)  $ (that is,
the permutation $w$ with $w\left(  i\right)  =j,w\left(  j\right)  =i$ and
$w\left(  k\right)  =k$ for $k\neq i,j$). The action of $S_{N}$ on coordinates
is defined by $\left(  xw\right)  _{i}=x_{w\left(  i\right)  }$ and is
extended to polynomials by $\left(  wp\right)  \left(  x\right)  :=p\left(
xw\right)  $ with the effect that $w\left(  x^{\alpha}\right)  =x^{w\alpha}$.
The operators $\mathcal{U}_{i}$ for $1\leq i\leq N$ are defined by
\[
\mathcal{U}_{i}p\left(  x\right)  :=\frac{\partial}{\partial x_{i}}\left(
x_{i}p\left(  x\right)  \right)  +\kappa\sum_{j=1,j\neq i}^{N}\frac
{x_{i}p\left(  x\right)  -x_{j}p\left(  x\left(  i,j\right)  \right)  }%
{x_{i}-x_{j}}-\kappa\sum_{j=1}^{i-1}p\left(  x\left(  i,j\right)  \right)  ,
\]
where $p$ is a polynomial ($\in\mathbb{Q}\left(  \kappa\right)  \left[
x_{1},x_{2},\ldots,x_{N}\right]  $). Then (see \cite[pp.291-2]{DX} for
details) $\mathcal{U}_{i}\mathcal{U}_{j}=\mathcal{U}_{j}\mathcal{U}_{i}$ for
$1\leq i,j\leq N$ and there is a crucial triangularity (in the sense of
matrices) property: $\mathcal{U}_{i}x^{\alpha}=\xi_{i}\left(  \alpha\right)
x^{\alpha}+q_{\alpha,i}\left(  x\right)  $ where $q_{\alpha,i}\left(
x\right)  $ is a sum of terms of the form $\pm\kappa x^{\beta}$ with certain
$\beta\vartriangleleft\alpha$ and
\begin{align*}
\xi_{i}\left(  \alpha\right)   &  :=\left(  N-r\left(  \alpha,i\right)
\right)  \kappa+\alpha_{i}+1\\
&  =\left(  \alpha_{i}-r\left(  \alpha,i\right)  \kappa\right)  +\left(
N\kappa+1\right)  ,
\end{align*}
for $1\leq i\leq N$ and $\alpha\in\mathbb{N}_{0}^{N}$.

The existence of the nonsymmetric Jack polynomials follows from a theorem of
elementary linear algebra. Suppose $\left\{  A^{\left(  i\right)  }:1\leq
i\leq N\right\}  $ is a collection of pairwise commuting lower triangular
$M\times M$ matrices over a field. If for each pair $\left(  i,j\right)
,1\leq i<j\leq M$ there is at least one matrix $A^{\left(  k\right)  }$ such
that $A_{ii}^{\left(  k\right)  }\neq A_{jj}^{\left(  k\right)  }$, then there
exists a unique set of $M$ linearly independent simultaneous (column)
eigenvectors for $\left\{  A^{\left(  i\right)  }\right\}  $ with each
eigenvector of the form $\left(  0,\ldots,0,1,\ldots\right)  ^{T}$.
Equivalently, there is a unique lower triangular unipotent matrix $V$ such
that $V^{-1}A^{\left(  i\right)  }V$ is diagonal for each $i$. Now apply this
result to the action of $\left\{  \mathcal{U}_{i}:1\leq i\leq N\right\}  $ on
the spaces of homogeneous polynomials with the standard basis $\left\{
x^{\alpha}:\alpha\in\mathbb{N}_{0}^{N},\left\vert \alpha\right\vert
=k\right\}  $, ordered by $\vartriangleright$, for $k\in\mathbb{N}_{0}$. It is
clear that $\alpha,\beta\in\mathbb{N}_{0}^{N}$ and $\alpha\neq\beta$ implies
that $\xi_{i}\left(  \alpha\right)  \neq\xi_{i}\left(  \beta\right)  $ for any
$i$ with $\alpha_{i}\neq\beta_{i}$ and generic $\kappa$. Thus for each
$\alpha\in\mathbb{N}_{0}^{N}$ there is a unique polynomial, called the
\textit{nonsymmetric Jack polynomial,}%
\[
\zeta_{\alpha}\left(  x\right)  =x^{\alpha}+\sum_{\beta\vartriangleleft\alpha
}A_{\beta\alpha}x^{\beta},
\]
with coefficients $A_{\beta\alpha}\in\mathbb{Q}\left(  \kappa\right)  $ such
that
\[
\mathcal{U}_{i}\zeta_{\alpha}=\xi_{i}\left(  \alpha\right)  \zeta_{\alpha
},\text{ for }1\leq i\leq N.
\]
The coefficients, as rational functions of $\kappa$, can have poles only at
certain negative rational numbers, which in turn are linked to the critical
pairs $\left(  \alpha,\beta\right)  $. This is a sketch of the argument (for a
detailed proof see \cite{SP} ): by the triangularity property there are
coefficients $B_{\beta\alpha}\in\mathbb{Q}\left(  \kappa\right)  $ such that
$x^{\alpha}=\zeta_{\alpha}+\sum_{\beta\vartriangleleft\alpha}B_{\beta\alpha
}\zeta_{\beta}$ for each $\alpha\in\mathbb{N}_{0}^{N}$. Extend the field
$\mathbb{Q}\left(  \kappa\right)  $ by adjoining a formal transcendental $\nu$
and consider the operator (on polynomials with coefficients in $\mathbb{Q}%
\left(  \kappa,\nu\right)  $)%
\[
\mathcal{T}_{\alpha}:=\prod_{\beta\vartriangleleft\alpha}\frac{\sum_{i=1}%
^{N}\nu^{i}\left(  \mathcal{U}_{i}-\xi_{i}\left(  \beta\right)  \right)
}{\sum_{i=1}^{N}\nu^{i}\left(  \xi_{i}\left(  \alpha\right)  -\xi_{i}\left(
\beta\right)  \right)  }.
\]
Now apply $\mathcal{T}_{\alpha}$ to the expression for $x^{\alpha}$ and
obtain:%
\[
\prod_{\beta\vartriangleleft\alpha}\frac{\sum_{i=1}^{N}\nu^{i}\left(
\mathcal{U}_{i}-\xi_{i}\left(  \beta\right)  \right)  }{\sum_{i=1}^{N}\nu
^{i}\left(  \xi_{i}\left(  \alpha\right)  -\xi_{i}\left(  \beta\right)
\right)  }x^{\alpha}=\zeta_{\alpha}+\sum_{\beta\vartriangleleft\alpha}%
B_{\beta\alpha}\mathcal{T}_{\alpha}\zeta_{\beta}=\zeta_{\alpha}.
\]
Since the variable $\nu$ does not appear in $\zeta_{\alpha}$ the denominators
of the coefficients $B_{\beta\alpha}$ must come from reducible factors of the
form
\begin{align*}
\sum_{i=1}^{N}\nu^{i}\left(  \xi_{i}\left(  \alpha\right)  -\xi_{i}\left(
\beta\right)  \right)   &  =\sum_{i=1}^{N}\nu^{i}\left(  \left(  r\left(
\beta,i\right)  -r\left(  \alpha,i\right)  \right)  \kappa+\alpha_{i}%
-\beta_{i}\right)  \\
&  =\left(  \sum_{i=1}^{N}c_{i}\nu^{i}\right)  \left(  m\kappa+n\right)
\end{align*}
with $c_{i}\in\mathbb{Q},~\beta\vartriangleleft\alpha$ and $m,n\geq1$. This is
exactly the property that $\left(  \alpha,\beta\right)  $ is a $\left(
-\frac{n}{m}\right)  $-critical pair.

By combinatorial means Knop and Sahi \cite{KS}\ showed that all coefficients
of $h\left(  \alpha,\kappa+1\right)  \zeta_{\alpha}^{x}$ are in $\mathbb{N}%
_{0}\left[  \kappa\right]  $ (polynomials in $\kappa$ with nonnegative integer
coefficients) and the coefficient of $x_{k+1}x_{k+2}\ldots x_{k+l}$ in
$\zeta_{\alpha}$ is $l!\kappa^{l}/h\left(  \alpha,\kappa+1\right)  $ for
$k=\ell\left(  \alpha\right)  $ and $l=\left\vert \alpha\right\vert $.

We conclude from the above discussion that for any node $\left(  i,j\right)
,1\leq j\leq\alpha_{i}$ with $h\left(  \alpha,\kappa+1;i,j\right)  =\left(
L\left(  \alpha;i,j\right)  +1\right)  \kappa+\alpha_{i}+1-j$ there must exist
at least one $\beta$ such that $\left(  \alpha,\beta\right)  $ is $-\left(
\alpha_{i}+1-j\right)  /\left(  L\left(  \alpha;i,j\right)  +1\right)
$-critical. The main reason for setting up the machinery of critical pairs is
to provide a tool for analyzing the dependence of the poles (as functions of
$\kappa$) in the coefficients of $\zeta_{\alpha}$ on the number of variables.
This will be illustrated in the last section.

\section{The Construction of Critical Pairs}

The main difficulty in extending the method from partitions to compositions is
to deal with tied values. Recall that for $\alpha_{i}=\alpha_{j}$ we have
$r\left(  \alpha,i\right)  <r\left(  \alpha,j\right)  $ if and only if $i<j$.
The definition of leg-length $L$ is more subtle for compositions. The
introduction of small deformations in the values makes it possible to use
essentially the same method as for partitions. Loosely speaking we use an
infinitesimal quantity $\upsilon$ which satisfies $0<i\upsilon<1$ for all
$i\geq1$, but of course the inequality will only be needed for all $i\leq N$
for some $N\geq\max\left(  \ell\left(  \alpha\right)  ,\ell\left(
\beta\right)  \right)  $. In the sequel we let $N=\ell\left(  \alpha\right)
+\left\vert \alpha\right\vert $ which suffices by Proposition \ref{lgt} and we
let $\upsilon:=\frac{1}{N+1}$.

\begin{definition}
For $\alpha\in\mathbb{N}_{0}^{N}$ let $\widetilde{\alpha}\in\mathbb{Q}^{N}$ be
given by $\widetilde{\alpha}_{i}:=\alpha_{i}-i\upsilon$ for $1\leq i\leq N$.
\end{definition}

\begin{proposition}
\label{a-tilde}For any $\alpha\in\mathbb{N}_{0}^{N}$ the following
hold:\newline(i) if $i\neq j$ then $\widetilde{\alpha}_{i}-\widetilde{\alpha
}_{j}\notin\mathbb{Z}$, in particular, $\widetilde{\alpha}_{i}\neq
\widetilde{\alpha}_{j}$,\newline(ii) $r\left(  \alpha,i\right)  =r\left(
\widetilde{\alpha},i\right)  $ for $1\leq i\leq N$,\newline(iii) there is a
unique permutation $w$ of $\left\{  1,\ldots,N\right\}  $ such that $r\left(
\alpha,w\left(  i\right)  \right)  =i$.
\end{proposition}

\begin{proof}
To show part (i) suppose $\alpha_{i}>\alpha_{j}$ then $\widetilde{\alpha}%
_{i}-\widetilde{\alpha}_{j}=\alpha_{i}-\alpha_{j}-\left(  i-j\right)
\upsilon\geq1-\left(  i-j\right)  \upsilon>0$, or suppose $\alpha_{i}%
=\alpha_{j}$ and $i<j$ then $\widetilde{\alpha}_{i}-\widetilde{\alpha}%
_{j}=\left(  j-i\right)  \upsilon>0$. In both cases, $\widetilde{\alpha}%
_{i}-\widetilde{\alpha}_{j}\notin\mathbb{Z}$. This shows that $r\left(
\widetilde{\alpha},i\right)  =\#\left\{  j:\widetilde{\alpha}_{j}%
>\widetilde{\alpha}_{i}\right\}  +1$ (extending the definition of the rank
function to elements of $\mathbb{Q}^{N}$). Now%
\begin{align*}
r\left(  \alpha,i\right)   &  =\#\left\{  j:j>i,\alpha_{j}>\alpha_{i}\right\}
+\#\left\{  j:j<i,\alpha_{j}\geq\alpha_{i}\right\}  +1\\
&  =\#\left\{  j:j>i,\widetilde{\alpha}_{j}>\widetilde{\alpha}_{i}\right\}
+\#\left\{  j:j<i,\widetilde{\alpha}_{j}>\widetilde{\alpha}_{i}\right\}  +1\\
&  =r\left(  \widetilde{\alpha},i\right)
\end{align*}
for $1\leq i\leq N$. The fact that the sorting permutation $w$ is unique
follows trivially from part (i).
\end{proof}

As noted before, part (iii) implies that $\alpha_{i}^{+}=\alpha_{w\left(
i\right)  }$, also that $\alpha_{i}=\alpha_{j}$ and $i<j$ implies $w\left(
i\right)  <w\left(  j\right)  $ (consider this as the formal proof that
$i\mapsto r\left(  \alpha,i\right)  $ is one-to-one). Here is the formula for
leg-length in terms of $\widetilde{\alpha}$.

\begin{proposition}
For $\alpha\in\mathbb{N}_{0}^{N}$ and for $1\leq i\leq\ell\left(
\alpha\right)  ,1\leq j\leq\alpha_{i}$ the leg-length satisfies the equation%
\[
L\left(  \alpha;i,j\right)  =\#\left\{  l:j-i\upsilon-1<\widetilde{\alpha}%
_{l}<\widetilde{\alpha}_{i}\right\}  .
\]

\end{proposition}

\begin{proof}
For $l<i$ the inequalities are $j-i\upsilon-1<\alpha_{l}-l\upsilon$ and
$\alpha_{l}-l\upsilon<\alpha_{i}-i\upsilon$, equivalent to $\alpha
_{l}+1-j>-\left(  i-l\right)  \upsilon>-1$ (hence $\alpha_{l}+1-j\geq0$) and
$\alpha_{i}-\alpha_{l}>\left(  i-l\right)  \upsilon>0$ (hence $\alpha
_{i}-\alpha_{l}\geq1$), respectively. For $l>i$ the inequalities $\alpha
_{l}+1-j>\left(  l-i\right)  \upsilon$ and $\alpha_{i}-\alpha_{l}>-\left(
l-i\right)  \upsilon$ are equivalent to $\alpha_{l}+1-j\geq1$ and $\alpha
_{i}-\alpha_{l}\geq0$ respectively.
\end{proof}

We begin the construction for a given $\alpha$ and the hook-length for a node
in the largest part. We use the permutation $w$ described in part (iii) of
Proposition \ref{a-tilde}. Suppose $1\leq n\leq\alpha_{w\left(  1\right)  }$
and $L\left(  \alpha;w\left(  1\right)  ,\alpha_{w\left(  1\right)
}+1-n\right)  =m-1$. We construct $\beta$ so that $\left(  \alpha
,\beta\right)  $ is a $\left(  -\frac{n}{m}\right)  $-critical pair.

There is another characterization of $m$: it has the property that
$\widetilde{\alpha}_{w\left(  m\right)  }>\widetilde{\alpha}_{w\left(
1\right)  }-n>\widetilde{\alpha}_{w\left(  j\right)  }$ for $j>m$; note%
\begin{align*}
m-1  &  =\#\left\{  l:\alpha_{w\left(  1\right)  }+1-n-w\left(  1\right)
\upsilon-1<\widetilde{\alpha}_{l}<\widetilde{\alpha}_{w\left(  1\right)
}\right\} \\
&  =\#\left\{  l:\widetilde{\alpha}_{w\left(  1\right)  }-n<\widetilde{\alpha
}_{l}<\widetilde{\alpha}_{w\left(  1\right)  }\right\}  .
\end{align*}
This shows that the $m$ largest parts of $\alpha$, excluding $\widetilde
{\alpha}_{w\left(  1\right)  }$, are in the latter interval; the equation
$\widetilde{\alpha}_{l}=\widetilde{\alpha}_{w\left(  1\right)  }-n$ is
impossible by part (i) of Proposition \ref{a-tilde}. We define an associated
sequence used in the construction: for $1\leq i\leq m$ and $k\geq0$ let%
\[
\xi_{mk+i}:=\widetilde{\alpha}_{w\left(  i\right)  }-nk.
\]

\begin{lemma}
The sequence $\left(  \xi_{i}\right)  _{i=1}^{N}$ is strictly decreasing and
$\xi_{j+1}<\widetilde{\alpha}_{w\left(  j\right)  }$ for large enough $j$.
There is a unique $T$ such that $\widetilde{\alpha}_{w\left(  m+s\right)
}<\xi_{m+s+1}$ for $1\leq s<T$ and $\widetilde{\alpha}_{w\left(  m+T\right)
}>\xi_{m+T+1}$ (the case of equality is ruled out).
\end{lemma}

\begin{proof}
The decreasing property has two cases. If $1\leq i<m$ then $\xi_{mk+i}%
-\xi_{mk+i+1}=\left(  \widetilde{\alpha}_{w\left(  i\right)  }-nk\right)
-\left(  \widetilde{\alpha}_{w\left(  i+1\right)  }-nk\right)  >0$. If $i=m$
then $\xi_{mk+m}-\xi_{mk+m+1}=\left(  \widetilde{\alpha}_{w\left(  m\right)
}-nk\right)  -\left(  \widetilde{\alpha}_{w\left(  1\right)  }-n\left(
k+1\right)  \right)  =\widetilde{\alpha}_{w\left(  m\right)  }-\left(
\widetilde{\alpha}_{w\left(  1\right)  }-n\right)  >0$. Since $\widetilde
{\alpha}_{w\left(  i\right)  }-\widetilde{\alpha}_{w\left(  j\right)  }%
\notin\mathbb{Z}$ for $i<j$ it is impossible for $\widetilde{\alpha}_{w\left(
m+s\right)  }=\xi_{m+s+1}$ when $s\geq1$.

Let $T_{0}=\sum_{i=1}^{m}\left\lfloor \alpha_{w\left(  i\right)
}/n\right\rfloor $. We claim $\xi_{m+T_{0}+1}<-1$. Set $\alpha_{w\left(
i\right)  }=nq_{i}+r_{i}$ with $q_{i}\geq0$ and $0\leq r_{i}\leq n-1$ for
$1\leq i\leq m$. Then $q_{1}\geq1$ since $n\leq\alpha_{w\left(  1\right)  }$
and $q_{1}\geq q_{i}\geq q_{1}-1$. This follows from the inequalities
$\alpha_{w\left(  1\right)  }\geq\alpha_{w\left(  i\right)  }\geq
\alpha_{w\left(  i\right)  }-n$, that is,%
\begin{align*}
nq_{1}+r_{1}  &  \geq nq_{i}+r_{i}\geq\left(  n-1\right)  q_{1}+r_{1},\\
r_{1}-r_{i}  &  \geq n\left(  q_{i}-q_{1}\right)  \geq r_{1}-r_{i}-n,
\end{align*}
but $1-n\leq r_{1}-r_{i}\leq n-1$ so $1-2n\leq n\left(  q_{i}-q_{1}\right)
\leq n-1$. Furthermore $q_{i}\geq q_{i+1}$ for $1\leq i<m$ since
$\alpha_{w\left(  i\right)  }\geq\alpha_{w\left(  i+1\right)  }$. Thus there
exists $k$ with $1<l\leq m$ so that $q_{i}=q_{1}$ for $1\leq i\leq l$ and
$q_{i}=q_{1}-1$ for $l<i\leq m$. Write $T_{0}+m+1=mk+i$ with $k\geq0$ and
$1\leq i\leq m$ . Then $mk+i=\left(  lq_{1}+\left(  m-l\right)  \left(
q_{1}-1\right)  \right)  +m+1=mq_{1}+l+1$ and so $i=\left(  l+1\right)
\operatorname{mod}m$. If $l<m$ then $k=q_{1},i=l+1$ and $\xi_{T_{0}%
+m+1}=\alpha_{w\left(  l+1\right)  }-nq_{1}-w\left(  l+1\right)
\upsilon=r_{l+1}-n-w\left(  l+1\right)  \upsilon<-1$. If $l=m$ then
$k=q_{1}+1,i=1$ and $\xi_{T_{0}+m+1}=\alpha_{w\left(  1\right)  }-\left(
n+1\right)  q_{1}-w\left(  1\right)  \upsilon=r_{1}-n-w\left(  1\right)
\upsilon<-1$. Finally $T_{0}\leq\frac{\left\vert \alpha\right\vert }{n}$ and
$m+T_{0}+1\leq\ell\left(  \alpha\right)  +\left\vert \alpha\right\vert +1$;
also $w\left(  j\right)  =j$ for $j>\ell\left(  \alpha\right)  $ and thus
$\widetilde{\alpha}_{w\left(  j\right)  }>-1>\xi_{j+1}$ for all sufficiently
large $j\leq N$. The existence and uniqueness of $T$ is now obvious. By part
(i) of Proposition \ref{a-tilde} $\widetilde{\alpha}_{w\left(  m+s\right)
}=\xi_{m+s+1}$ is impossible for $s\geq1$.
\end{proof}

\begin{definition}
\label{defbeta}With $\alpha$ and $T$ as described above let $t:=\left(
\left(  T-1\right)  \operatorname{mod}m\right)  +1$, $k:=\frac{T-t}{m}$ (so
that $T=mk+t$ and $1\leq t\leq m$) and define $\beta\in\mathbb{N}_{0}^{N}$ by%
\begin{equation}
\beta_{w\left(  i\right)  }:=\left\{
\begin{tabular}
[c]{ll}%
$\alpha_{w\left(  i\right)  }-\left(  k+1\right)  n,$ & $1\leq i\leq t$\\
$\alpha_{w\left(  i\right)  }-kn,$ & $t<i\leq m$\\
$\alpha_{w\left(  i\right)  }+n,$ & $m+1\leq i\leq m+T$\\
$\alpha_{w\left(  i\right)  },$ & $m+T<i.$%
\end{tabular}
\right.  \label{eq-beta}%
\end{equation}

\end{definition}

The following is the main result. The notations $\alpha,\beta,\xi,w,m,n,T,t$
continue with the definitions given above. The proof is broken up in several
lemmas. It is possible that $m=1$, in which case one makes the obvious
modifications in the following statements.

\begin{theorem}
For $\alpha\in\mathbb{N}_{0}^{N}$ the composition $\beta$ in Definition
\ref{defbeta} has the property that $m\kappa+n$ divides $\left(  r\left(
\beta,w\left(  i\right)  \right)  -i\right)  \kappa+\alpha_{w\left(  i\right)
}-\beta_{w\left(  i\right)  }$ for all $i$, and $\alpha\vartriangleright\beta
$, that is, $\left(  \alpha,\beta\right)  $ is a $\left(  -\frac{n}{m}\right)
$-critical pair.
\end{theorem}

Before we present the details of the proof we explain how the Theorem can be
used as an algorithm. Here is an informal description.

\begin{algorithm}
Start with $\alpha,m,n,N=\left\vert \alpha\right\vert +\ell\left(
\alpha\right)  ,\upsilon$ as described above.

\begin{enumerate}
\item Compute the permutation $w\in S_{N}$ with the property $r\left(
\alpha,w\left(  i\right)  \right)  =i$ for $1\leq i\leq N$,

\item for $i=1,2,\ldots,m$ set $\xi_{i}:=\widetilde{\alpha}_{w\left(
i\right)  }$, and set $\xi_{m+1}:=\xi_{1}-n$,

\item for $s=2,3,\ldots$ set $\xi_{m+s}:=\xi_{s}-n$ until $\xi_{m+s}%
<\widetilde{\alpha}_{w\left(  m+s-1\right)  }$,

\item set $T:=s-1$ (where $s$ is the first value in step 3 for which
$\xi_{m+s}<\widetilde{\alpha}_{w\left(  m+s-1\right)  }$),

\item set $t:=\left(  \left(  T-1\right)  \operatorname{mod}m\right)  +1$,
$k:=\dfrac{T-t}{m}$,

\item use Equation \ref{eq-beta} to compute $\beta$.
\end{enumerate}
\end{algorithm}

Here is an example: let $\alpha=\left(  0,3,5,6,6,1\right)  $; then $h\left(
\alpha,\kappa+1;4,4\right)  =4\kappa+3$ and $\left(  w\left(  i\right)
\right)  _{i=1}^{N}=\left(  4,5,3,2,6,1,7,8,\ldots\right)  $. The sequence
$\left(  \widetilde{\alpha}_{w\left(  i\right)  }\right)  _{i=1}^{N}$ is%
\[
\left(  6-4\upsilon,6-5\upsilon,5-3\upsilon,3-2\upsilon,1-6\upsilon
,-\upsilon,-7\upsilon,-8\upsilon,-9\upsilon,-10\upsilon,-11\upsilon
,\ldots\right)
\]
and the sequence $\left(  \xi_{i}\right)  _{i=1}^{N}$ is computed up to the
$11^{th}$ term
\[
\left(  6-4\upsilon,6-5\upsilon,5-3\upsilon,3-2\upsilon,3-4\upsilon
,3-5\upsilon,2-3\upsilon,-2\upsilon,-4\upsilon,-5\upsilon,-1-3\upsilon
,\ldots\right)
\]
since $\widetilde{\alpha}_{w\left(  4+6\right)  }>\xi_{4+7}$ and
$\widetilde{\alpha}_{w\left(  4+s\right)  }<\xi_{5+s}$ for 1$\leq s<6$. So
$T=6,t=2,k=1$ and $\left(  \beta_{w\left(  i\right)  }\right)  _{i=1}%
^{N}=\left(  0,0,2,0,4,3,3,3,3,3,0,\ldots\right)  $, $\beta=\left(
3,0,2,0,0,4,3,3,3,3\right)  $. Now $\left(  r\left(  \beta,i\right)  \right)
_{i=1}^{10}=\left(  2,8,7,9,10,1,3,4,5,6\right)  $ and $\left(  \alpha
,\beta\right)  $ is indeed a $\left(  -\frac{3}{4}\right)  $-critical pair.

It is possible that $T=1$, for example let $\alpha=\left(  2,6,5,2\right)  $
then $h\left(  \alpha,\kappa+1;2,4\right)  =2\kappa+3,\left(  \widetilde
{\alpha}_{w\left(  i\right)  }\right)  _{i=1}^{4}=\left(  6-2\upsilon
,5-3\upsilon,2-\upsilon,2-4\upsilon\right)  $ and $\xi=\left(  6-2\upsilon
,5-3\upsilon,3-2\upsilon,2-3\upsilon\right)  $ so that $\widetilde{\alpha
}_{w\left(  2+1\right)  }>\xi_{4}$; then $\beta=\left(  5,3,5,2\right)  $.

The proof of the theorem is broken up into several lemmas.

\begin{lemma}
The following inequalities hold:\newline(i) $\widetilde{\alpha}_{w\left(
m+T+1\right)  }<\xi_{T+m}<\xi_{T+1}<\widetilde{\alpha}_{w\left(  m+T\right)
}+n$ (omit $\xi_{T+m}$ if $m=1$)\newline(ii) for $t<m,$ $\widetilde{\beta
}_{w\left(  m+T\right)  }>\widetilde{\beta}_{w\left(  t+1\right)  }%
>\ldots>\widetilde{\beta}_{w\left(  m\right)  }>\widetilde{\beta}_{w\left(
1\right)  }>\ldots>\widetilde{\beta}_{w\left(  t\right)  }>\widetilde{\beta
}_{w\left(  m+T+1\right)  }$\newline(iii) for $t=m,$ $\widetilde{\beta
}_{w\left(  m+T\right)  }>\widetilde{\beta}_{w\left(  1\right)  }%
>\ldots>\widetilde{\beta}_{w\left(  m\right)  }>\widetilde{\beta}_{w\left(
m+T+1\right)  }$.
\end{lemma}

\begin{proof}
By construction $\widetilde{\alpha}_{w\left(  m+T\right)  }+n>\xi
_{m+T+1}+n=\left(  \xi_{T+1}-n\right)  +n>\xi_{T+m}$, by the decreasing
property of $\left\{  \xi_{j}\right\}  $. If $T>1$ then $\xi_{T+m}%
>\widetilde{\alpha}_{w\left(  m+T-1\right)  }>\widetilde{\alpha}_{w\left(
m+T+1\right)  }$. If $T=1$ then $\widetilde{\alpha}_{w\left(  m+1\right)
}+n>\xi_{2}>\xi_{m+1}>\widetilde{\alpha}_{w\left(  m+1\right)  }%
>\widetilde{\alpha}_{w\left(  m+2\right)  }$. This proves part (i). As before
$T=km+t$. By construction,%
\[
\widetilde{\beta}_{w\left(  i\right)  }=\left\{
\begin{tabular}
[c]{ll}%
$\xi_{mk+m+i}$, & $1\leq i\leq t$\\
$\xi_{mk+i},$ & $t<i\leq m$\\
$\widetilde{\alpha}_{w\left(  i\right)  }+n,$ & $m+1\leq i\leq m+T$\\
$\widetilde{\alpha}_{w\left(  i\right)  },$ & $m+T<i.$%
\end{tabular}
\right.  .
\]
Thus the inequality in part (i) shows $\widetilde{\beta}_{w\left(
m+T+1\right)  }<\widetilde{\beta}_{w\left(  t\right)  }<\widetilde{\beta
}_{w\left(  t+1\right)  }<\widetilde{\beta}_{w\left(  m+T\right)  }$.
\end{proof}

\begin{lemma}
The following rank values hold:\newline(i) $r\left(  \beta,w\left(
m+i\right)  \right)  =i$ for $1\leq i\leq T,$\newline(ii) $r\left(
\beta,w\left(  t+i\right)  \right)  =mk+t+i$ for $1\leq i\leq m-t,$%
\newline(iii) $r\left(  \beta,w\left(  i\right)  \right)  =m\left(
k+1\right)  +i$ for $1\leq i\leq t,$\newline(iv) $r\left(  \beta,w\left(
i\right)  \right)  =i$ for $m+T+1\leq i\leq N$.
\end{lemma}

\begin{proof}
By parts (ii) and (iii) of the previous lemma, in decreasing order the $m$
largest values of $\widetilde{\beta}$ are $\widetilde{\beta}_{w\left(
m+1\right)  },\ldots,\widetilde{\beta}_{w\left(  m+T\right)  }$, the next
$m-t$ values are $\widetilde{\beta}_{w\left(  t+1\right)  },\ldots
,\widetilde{\beta}_{w\left(  m\right)  }$, the next $t$ values are
$\widetilde{\beta}_{w\left(  1\right)  },\ldots,\widetilde{\beta}_{w\left(
t\right)  }$, and the remaining are $\widetilde{\beta}_{w\left(  m+T+1\right)
},\ldots$. Note that $r\left(  \beta,w\left(  t+i\right)  \right)
=T+i=mk+t+i$ in part (ii) and $r\left(  \beta,w\left(  i\right)  \right)
=T+\left(  m-t\right)  +i=m\left(  k+1\right)  +i$.
\end{proof}

\begin{lemma}
The composition $\beta$ satisfies the condition that $m\kappa+n$ divides
\newline$\left(  r\left(  \beta,w\left(  i\right)  \right)  -i\right)
\kappa+\alpha_{w\left(  i\right)  }-\beta_{w\left(  i\right)  }$ for all
$i\leq N$.
\end{lemma}

\begin{proof}
Let $\gamma_{i}=\left(  r\left(  \beta,w\left(  i\right)  \right)  -i\right)
\kappa+\alpha_{w\left(  i\right)  }-\beta_{w\left(  i\right)  }$. For $1\leq
i\leq t$, $\gamma_{i}=\left(  m\left(  k+1\right)  +i-i\right)  \kappa
+\alpha_{w\left(  i\right)  }-\left(  \alpha_{w\left(  i\right)  }-\left(
k+1\right)  n\right)  =\left(  k+1\right)  \left(  m\kappa+n\right)  $. For
$t<i\leq m$, $\gamma_{i}=mk\kappa+\alpha_{w\left(  i\right)  }-\left(
\alpha_{w\left(  i\right)  }-kn\right)  =k\left(  m\kappa+n\right)  $. For
$m+1\leq i\leq m+T$, $\gamma_{i}=\left(  \left(  i-m\right)  -i\right)
\kappa+\alpha_{w\left(  i\right)  }-\left(  \alpha_{w\left(  i\right)
}+n\right)  =-\left(  m\kappa+n\right)  $, and $\gamma_{i}=0$ for $i>m+T$.
\end{proof}

We must show that $\alpha\vartriangleright\beta$ to complete the proof of the
theorem. For $1\leq i\leq N$ let $\varepsilon\left(  i\right)  \in
\mathbb{N}_{0}^{N}$ denote the standard basis element, that is, $\varepsilon
\left(  i\right)  _{j}=\delta_{ij}$. The idea is to describe the construction
as a sequence of compositions, each of which is produced by adding $n\left(
\varepsilon\left(  w\left(  m+i\right)  \right)  -\varepsilon\left(  w\left(
\left(  i-1\right)  \operatorname{mod}m+1\right)  \right)  \right)  $ to the
previous one, for $i=1,2,\ldots T$. The effect of the argument $\left(
\left(  i-1\right)  \operatorname{mod}m\right)  +1$ is to cycle through the
values $w\left(  1\right)  ,\ldots,w\left(  m\right)  $.

\begin{lemma}
Let $\beta^{\left(  s\right)  }:=\alpha-n\sum_{i=1}^{s}\left(  \varepsilon
\left(  w\left(  \left(  i-1\right)  \operatorname{mod}m+1\right)  \right)
-\varepsilon\left(  w\left(  m+i\right)  \right)  \right)  $ for $0\leq s\leq
T$. Then $\beta^{\left(  0\right)  }=\alpha,\beta^{\left(  T\right)  }=\beta$
and $\beta^{\left(  s\right)  }\vartriangleright\beta^{\left(  s+1\right)  }$
for $0\leq s<T$.
\end{lemma}

\begin{proof}
We use Lemma 8.2.3 from \cite[p.289]{DX}. This states that if $\lambda$ is a
partition such that $1\leq n<\lambda_{i}-\lambda_{j}$ (for some $i<j$) then
$\lambda\succ\left(  \lambda-n\left(  \varepsilon\left(  i\right)
-\varepsilon\left(  j\right)  \right)  \right)  ^{+}$. As a consequence we
have that if $\gamma\in\mathbb{N}_{0}^{N}$ and $1\leq n<\gamma_{i}-\gamma_{j}$
for some $i,j$ then $\gamma^{+}\succ\left(  \gamma-n\left(  \varepsilon\left(
i\right)  -\varepsilon\left(  j\right)  \right)  \right)  ^{+}$, that is,
$\gamma\vartriangleright\left(  \gamma-n\left(  \varepsilon\left(  i\right)
-\varepsilon\left(  j\right)  \right)  \right)  $. Suppose $s=ml+i<m+T$ with
$0\leq i<m$ and $l\geq0$. Then $\beta_{w\left(  j\right)  }^{\left(  s\right)
}=\alpha_{w\left(  j\right)  }+n$ for $m+1\leq j\leq m+s$, $\beta_{w\left(
j\right)  }^{\left(  s\right)  }=\alpha_{w\left(  j\right)  }-n\left(
l+1\right)  $ for $1\leq j\leq i$, $\beta_{w\left(  j\right)  }^{\left(
s\right)  }=\alpha_{w\left(  j\right)  }-nl$ for $i+1\leq j\leq m$ and
$\beta_{w\left(  m+s+1\right)  }=\alpha_{w\left(  m+s+1\right)  }$. By
definition, $\beta^{\left(  s+1\right)  }=\beta^{\left(  s\right)  }-n\left(
\varepsilon\left(  w\left(  i+1\right)  \right)  -\varepsilon\left(  w\left(
m+s+1\right)  \right)  \right)  $. Let $\delta_{s}$ denote the difference
between the two affected values, that is,
\begin{align*}
\delta_{s}  &  :=\beta_{w\left(  i+1\right)  }^{\left(  s\right)  }%
-\beta_{w\left(  m+s+1\right)  }^{\left(  s\right)  }\\
&  =\alpha_{w\left(  i+1\right)  }-nl-\alpha_{w\left(  m+s+1\right)  }\\
&  =\xi_{ml+i+1}-\widetilde{\alpha}_{w\left(  m+s+1\right)  }+\left(  w\left(
s+1\right)  -w\left(  m+s+1\right)  \right)  \upsilon.
\end{align*}
By the construction of the Algorithm $\xi_{m+s+1}>\widetilde{\alpha}_{w\left(
m+s\right)  }$ for $s<T$ and $\widetilde{\alpha}_{w\left(  m+s\right)
}>\widetilde{\alpha}_{w\left(  m+s+1\right)  }$. Thus $\delta_{s}%
-n=\xi_{m+s+1}-\widetilde{\alpha}_{w\left(  m+s+1\right)  }+\left(  w\left(
s+1\right)  -w\left(  m+s+1\right)  \right)  \upsilon>\left(  w\left(
s+1\right)  -w\left(  m+s+1\right)  \right)  \upsilon>-1$. If $\delta
_{s}-n\geq1$ then by the above argument we have $\beta^{\left(  s\right)
}\vartriangleright\beta^{\left(  s+1\right)  }$. If $\delta_{s}-n=0$ then
$\beta_{w\left(  m+s+1\right)  }^{\left(  s+1\right)  }=\beta_{w\left(
i+1\right)  }^{\left(  s\right)  },\beta_{w\left(  i+1\right)  }^{\left(
s+1\right)  }=\beta_{w\left(  m+s+1\right)  }^{\left(  s\right)  }$ and
$w\left(  s+1\right)  <w\left(  m+s+1\right)  $; thus $\left(  \beta^{\left(
s+1\right)  }\right)  ^{+}=\left(  \beta^{\left(  s\right)  }\right)  ^{+}$
and $\beta^{\left(  s\right)  }\succ\beta^{\left(  s+1\right)  }$ (for any
composition $\gamma$, if $i<j$ and $\gamma_{i}>\gamma_{j}$ then $\gamma
\succ\left(  i,j\right)  \gamma$), that is, $\beta^{\left(  s\right)
}\vartriangleright\beta^{\left(  s+1\right)  }$.
\end{proof}

This completes the proof of the theorem when the hook length $m\kappa+n$ is
associated with the largest part (either $\alpha_{w\left(  1\right)  }%
>\alpha_{w\left(  2\right)  }$ or $\alpha_{w\left(  1\right)  }=\alpha
_{w\left(  i\right)  }$ and $i\neq1$ implies $w\left(  1\right)  <w\left(
i\right)  $). Suppose for some $l\geq1$ that $L\left(  \alpha;w\left(
l+1\right)  ,\alpha_{w\left(  l+1\right)  }+1-n\right)  =m-1$. Then the
algorithm is applied with the arguments of $w$ and the ranks all shifted by
$l$; and the largest $l$ parts are not changed. Thus $\beta_{w\left(
i\right)  }:=\alpha_{w\left(  i\right)  }$ for $1\leq i\leq l$ and
\[
\beta_{w\left(  l+i\right)  }:=\left\{
\begin{tabular}
[c]{ll}%
$\alpha_{w\left(  l+i\right)  }-\left(  k+1\right)  n,$ & $1\leq i\leq t$\\
$\alpha_{w\left(  l+i\right)  }-kn,$ & $t<i\leq m$\\
$\alpha_{w\left(  l+i\right)  }+n,$ & $m+1\leq i\leq m+T$\\
$\alpha_{w\left(  l+i\right)  },$ & $m+T<i,$%
\end{tabular}
\ \right.
\]
using the same notations $k,T,t$ as above.

\section{Examples and Discussion}

The first example is a partition type: $\alpha=\left(
9,8,8,7,4,3,3,2,2\right)  ,n=3,m=4,T=9$ and $\beta=\left(
0,2,2,1,7,6,6,5,5,3,3,3,3\right)  $. The ranks of $\beta$ are $\left(
r\left(  \beta,i\right)  \right)  _{i=1}^{13}=\left(
13,10,11,12,1,2,3,4,5,6,7,8,9\right)  $. To illustrate the situation where the
node $\left(  i,j\right)  $ (with hook length $h\left(  \alpha,\kappa
+1;i,j\right)  =\left(  L\left(  \alpha;i,j\right)  +1\right)  \kappa
+\alpha_{i}+1-j$) is not in the row of rank $1$, consider the hook at $\left(
2,5\right)  $ in $\alpha=\left(  9,8,8,5,4,4\right)  $, then $L\left(
\alpha;2,5\right)  =2,n=4$ and $\beta=\left(  9,4,4,5,8,8\right)  $ with ranks
$\left(  1,5,6,4,2,3\right)  $. The next example is a composition
$\alpha=\left(  0,3,5,6,6,4,1\right)  $ (with $\left(  r\left(  \alpha
,i\right)  \right)  _{i=1}^{7}=\left(  7,5,3,1,2,4,6\right)  $) with $n=3$;
then $m=L\left(  \alpha;4,4\right)  +1=5,T=6$ and $\beta=\left(
3,0,2,0,0,1,4,3,3,3,3,3\right)  $ (with $\left(  r\left(  \beta,i\right)
\right)  _{i=1}^{12}=\left(  2,10,8,11,12,9,1,3,4,5,6,7\right)  $).

There is an analogous situation for critical pairs when $\beta$ is the given
composition; in this case the expansion of nonsymmetric Jack polynomials in
terms of the $p$-basis (see \cite[p.298]{DX}) suggests that the linear factors
$\left(  m\kappa+n\right)  $ of $h\left(  \beta,1\right)  $ lead to $\left(
-\frac{n}{m}\right)  $-critical pairs $\left(  \alpha,\beta\right)  $. However
we consider the uniqueness problem described below as more important for
applications, and will not further investigate $h\left(  \beta,1\right)  $.

A natural question occurs: for a given $\alpha$ and hook-length $m\kappa+n$
does one step of the algorithm produce all possible solutions for $\beta$? It
makes sense to consider a sequence of steps because there is a transitive
property for critical pairs: if $\left(  \alpha,\beta^{\left(  1\right)
}\right)  $ and $\left(  \beta^{\left(  1\right)  },\beta^{\left(  2\right)
}\right)  $ are $\left(  -\frac{n}{m}\right)  $-critical pairs then so is
$\left(  \alpha,\beta^{\left(  2\right)  }\right)  $ ($\vartriangleright$
being a partial order). If one step sufficed then there would be a uniqueness
result of the form: suppose the multiplicity of the linear factor $\left(
m\kappa+n\right)  $ in the hook-length product $h\left(  \alpha,\kappa
+1\right)  $ is one, then there is a unique $\beta$ so that $\left(
\alpha,\beta\right)  $ is a $\left(  -\frac{n}{m}\right)  $-critical pair.
However, this may fail if $m,n$ are not relatively prime: for $\alpha=\left(
6,3,1,1\right)  ,n=6,m=4$ the algorithm produces $\beta^{\left(  1\right)
}=\left(  0,3,1,1,6\right)  ;$the multiplicity of $\left(  2\kappa+3\right)  $
is one in both $h\left(  \alpha,\kappa+1\right)  $ and $h\left(
\beta^{\left(  1\right)  },\kappa+1\right)  $. Apply the algorithm to
$\beta^{\left(  1\right)  }$ with $n=3,m=2$ to obtain $\beta^{\left(
2\right)  }=\left(  0,3,4,1,3\right)  $. The process stops since $h\left(
\beta^{\left(  2\right)  },\kappa+1\right)  $ does not have $\left(
2\kappa+3\right)  $ as a factor. We conjecture there is uniqueness if $m,n$
are relatively prime (and the multiplicity is one); a particular case of this
was established in \cite{SP}. This is crucial because it is used to show that
$\zeta_{\alpha}$ has no pole at $\kappa=-\frac{n}{m}$ when the number of
variables is less than $\ell\left(  \beta\right)  $. Here is an example:
$\alpha=\left(  7,6,6,4,4\right)  $ with $n=2$, so that $m=3$, then the unique
$\beta=\left(  1,0,0,6,6,2,2,2,2,2,2,2\right)  $; the application is to
$\zeta_{\alpha}$ on $\mathbb{R}^{10}$ and the uniqueness of $\beta$ implies
that $\kappa=-\frac{2}{3}$ is not a pole of $\zeta_{\alpha}$ (for less than 12
variables). However the uniqueness result uses specific properties of a class
of partitions and $m,n$ are relatively prime. There is a weak uniqueness
result concerning the sign changes in the sequence $\left(  \widetilde{\alpha
}_{w\left(  m+i\right)  }-\xi_{m+i+1}\right)  _{i=1}^{N}$. Recall that $T$ is
chosen so that $\widetilde{\alpha}_{w\left(  m+s\right)  }-\xi_{m+s+1}<0$ for
$1\leq s<T$ and $\widetilde{\alpha}_{w\left(  m+T\right)  }-\xi_{m+T+1}>0$.

\begin{proposition}
Suppose for some $s>T$ that $\widetilde{\alpha}_{w\left(  m+s\right)  }%
-\xi_{m+s+1}>0$ and $\widetilde{\alpha}_{w\left(  m+s+1\right)  }-\xi
_{m+s+2}<0$, then the hook-length at the node $\left(  w\left(  i\right)
,\alpha_{w\left(  i\right)  }+1-nl\right)  $ is $l\left(  m\kappa+n\right)  $
where $m+s+1=ml+i$ and $1\leq i\leq m$.
\end{proposition}

\begin{proof}
By hypothesis $\widetilde{\alpha}_{w\left(  ml+i+1\right)  }>\widetilde
{\alpha}_{w\left(  i\right)  }-nl>\xi_{m+s+2}>\widetilde{\alpha}_{w\left(
ml+i\right)  }$. This implies $L\left(  \alpha;w\left(  i\right)
,\alpha_{w\left(  i\right)  }+1-nl\right)  =ml-1$.
\end{proof}

The Proposition implies that if there is only one hook-length divisible by
$m\kappa+n$ in the rows $w\left(  i\right)  $ for $1\leq i\leq m$ then
$\widetilde{\alpha}_{w\left(  m+i\right)  }>\xi_{m+i+1}$ for all $i\geq T$; so
there is only one sign-change.

It appears that there can be a considerably larger number of solutions than
the multiplicity. Here is an example: $\alpha=\left(  9,7,6,5,2\right)  $, the
multiplicity of $\left(  2\kappa+3\right)  $ in $h\left(  \alpha
,\kappa+1\right)  $ is 4. The relevant hook lengths are $2\kappa+3$ at nodes
$\left(  1,7\right)  $ and $\left(  3,4\right)  $, and $4\kappa+6$ at nodes
$\left(  1,4\right)  $ and $\left(  2,2\right)  $. The algorithm produces
$\left(  6,7,9,5,2\right)  ,\allowbreak\left(  9,7,0,2,5,3,3\right)
,\allowbreak\left(  3,7,6,5,8\right)  ,\allowbreak\left(
9,1,0,5,2,6,6\right)  $ respectively for these nodes. But one can continue the
process: for example the multiplicity of $\left(  2\kappa+3\right)  $ in
\newline$h\left(  \left(  6,7,9,5,2\right)  ,\kappa+1\right)  $ is 3, and the
algorithm (applied to $\left(  6,7,9,5,2\right)  $) produces three more
solutions for $\beta$, one being $\left(  6,7,3,5,8\right)  $. One could
speculate that there is a lattice of solutions, ordered by $\vartriangleright
$. Finally, one can ask if there are $\left(  -\frac{n}{m}\right)  $-critical
pairs without a corresponding hook-length $m\kappa+n$. It seems doubtful, but
we will leave this unanswered.

The author thanks the referees for useful comments leading to an improved presentation.


\begin{thebibliography}{9}                                                                                                %


\bibitem {SP}C. Dunkl, Singular polynomials for the symmetric groups,
\emph{Int. Math. Research Not.} \textbf{2004} (2004), \#67, 3607-3635, arXiv:math.RT/0403277.

\bibitem {DX}C. Dunkl and Y. Xu, \textit{Orthogonal Polynomials of Several
Variables}, Encycl. of Math. and its Applications \textbf{81}, Cambridge
University Press, Cambridge, 2001.

\bibitem {KS}F. Knop and S. Sahi, A recursion and a combinatorial formula for
Jack polynomials. \emph{Invent. Math.} \textbf{128} (1997), 9--22.
\end{thebibliography}
\end{document}